# SCIENTIFIC CALCULATOR WITH THE AID OF GEOMETRY AND BASED UPON IT, A MECHANICAL CALCULATOR

**Narinder K Wadhawan**[1],

[1]Civil Servant, Indian Administrative Service Retired,
House No. 563, Sector 2, Panchkula-134112, Haryana, India, Cellphone: +919888511248,
e-mails: narinderkw@gmail.com

Abstract: Scientific calculations involving operations of multiplication, division, exponents, inverse of exponents of real numbers, geometric mean, reciprocal, Euler's number, logarithm, antilogarithm of a real number, are generally carried out with battery operated electronic calculators. In this paper, geometric methods employing properties of similar right angled triangles have been innovated for performing error free scientific calculations without the use of battery and also based upon it, a mechanical analogue calculator has also been devised. A right angled triangle with perpendicular side of unit length, is drawn and from the point of where two sides subscribe a right angle, a perpendicular is drawn upon hypotenuse. From the point of intersection of this perpendicular and hypotenuse, another perpendicular is drawn upon base and this process is continued till '$n$' perpendiculars are drawn, the resultant figure so formed, is a right angled triangle with '$n$' similar triangles inscribed in it. Using properties of similar triangles, it is found, $1^{st}$, $2^{nd}$, $3^{rd}$, …, $n^{th}$, perpendiculars are in geometric series. If the length of $n^{th}$ perpendicular in such a triangle is equated with the given real positive quantity that is less than one and on adjusting base angle, length of first perpendicular will measure $n^{th}$ root of the given quantity. If the given quantity is more than one, its reciprocal is equated with the length of $n^{th}$ perpendicular, then reciprocal of length of first perpendicular will be $n^{th}$ root. In this way, exploiting geometric properties, others scientific calculations have been performed.

Keywords: Antilogarithm; Geometric Computation; Logarithm; Mathematical Analog Calculator; Mathematical Operations; Right Angled Triangle.
MSC Subject Classification: 68U05, 65D18.

## 1. INTRODUCTION

This research paper deals with inventing *geometrical methods and devising a mechanical analogue calculator* for performing operations of *multiplication*, *division of two real numbers except zero and infinity*, *geometric mean of two real positive numbers except zero and infinity, reciprocal of a real number except zero and infinity*, *computations of logarithm of a real positive number, antilogarithm of a real number except zero and infinity, computation of values of $x^n$, $x^{1/n}$, $x^{m/n}$, where $x, m,$ and $n$ are real numbers except zero or infinity, computation of values of $n, 1/n, m/n,$ when given values of $x^n$, $x^{1/n}$, $x^{m/n}$, are real positive numbers, and computation of value of Euler's number $e$*. For achieving these purposes, the characteristic properties of similar right angled triangles formed on dropping perpendiculars upon hypotenuse from various points on its base and also dropping perpendiculars from various points on its hypotenuse upon base as will be explained as this paper proceeds , were utilised.

Going down the history lane, *Abacus* was the first tool invented in *Sumeria* around 2500 B.C. for performing mathematical calculations (Martin 1992; Williams 1997). With the passage of time and advancement of technology, pocket electronic calculator was invented in the year 1970. To give background of analogue computation with electrical components vis a vis geometric computations with regular figures, some examples are forth hereinafter.

When a capacitor of value $C$ is connected across a voltage source $v_s$, charge starts storing in it till voltage across it say $v_c$, reaches the value of applied voltage $v_s$. If $i_t$ is the electric current following in it at any instant $t$, voltage across capacitor at that instant is given by the relation, $v_c = \int (i_t/C) \, dt$ (Roberge 2013). If maximum value of current $i_m$ is taken care of by use of resistance say $r$ connected in series with capacitor $C$, output voltage $v_c$ across capacitor is integration of current $i_t$. Conversely, if out put voltage say $v_r$ is taken across resistance $r$, $v_r$ is indicative of differentiation of voltage across capacitor $v_c$ (Roberge 2013). Also if an inductor of value $L$ and resistance of value $r$ in series are connected across a voltage source, then $v_L = -L(di/dt)$, where $v_L$ and $di/dt$ are voltage across inductor and rate of chage of current with time respectively (Roberge 2013). However, if only a resistance $R$ is connected across a voltage source, then voltage across resistance, $v_R = v_s = Ri$ or $i = v_s/R$ and that means current, in this case, is indicative of division of voltage $v_s$ by resistance $R$. If two or more resistances say $R_1, R_2, R_3, \ldots, R_n$ are connected in parallel across a voltage source $v_s$, then total current being drawn from the source is $v_s \sum_{i=1}^{n} 1/R_i$, where $i$ varies from 1 to $n$ is a classic case of summation of harmonic series provided source voltage is unity. If a diode and a resistance in series are connected across a voltage source that forward-biases the diode, the voltage across the diode is given by relation $V = K \ln (I/I_0)$. If physical quantities like temperature etcetera remain constant, then $K$ and $I_0$ can be assumed constant and in that case, voltage across diode is $V = K \ln (I/I_0)$, which is indicative of natural logarithm of current $I$ flowing through it and voltage across resistance is indicative of antilogarithm of voltage (Roberge 2013). Some transducers like piezo electric crystals bellows, bourdon tube, thermo couples respond to parameters vibrations, pressure, temperature, flow, vibrations likewise and the characteristics of these elements are used for measurement and control in *Instrumentation* (Morris and Langari 2021).

In the backdrop of above explanation and referring to Figure 1, mathematical pattern of variation in the length of perpendiculars, $AB, BD, DE, EF$ ..., drawn in triangle $ABC$ with right angle at $B$, base $BC$ and hypotenuse $AC$ are derived for performing error free calculations. These geometric properties were also used to calculate values of real root(s) of a quintic polynomial equations (Wadhawan 2023). While explaining logarithmic curves, this property of right angled triangles was used by Rene Descartes (Smith and Latham translation 1954). $\triangle ABC$ has its vertex at point $A$, base $BC$, $\angle ABC = 90$ degrees and $\angle ACB$ abbreviated $\angle C$ at point $C$. $\angle C$ can have any value less than 90 degrees and more than 0 degree. For *construction of angle C* corresponding to $\cos C = x$, where magnitude of $x$ is less than one and more than zero, a horizontal straight line $BC$ of length $x$ unit, will be drawn and a perpendicular $BA$, from point $B$ upon base $BC$ will be raised and from point $C$, an arc of unit radius will be drawn intersecting the perpendicular at point $A$, then angle $ACB$ or simply angle $C$ will correspond to $\cos C =$

$x$. Continuing with the construction, from point $B$, a $\perp BD$ is drawn on hypotenuse $AC$, from point $D$, a $\perp DE$ is drawn upon base $BC$ and so on. Lengths $BD, EF, GH, IJ, KL, \ldots$, given notations as $p_1, p_3, p_5, \ldots$, in the paper, are perpendiculars upon hypotenuse, $AC. DE, FG. HI, JK, LM, \ldots$, mentioned as $p_2, p_4, p_6, \ldots$, are perpendiculars upon base $BC$. This construction henceforth will not be repeated to avoid redundancy.

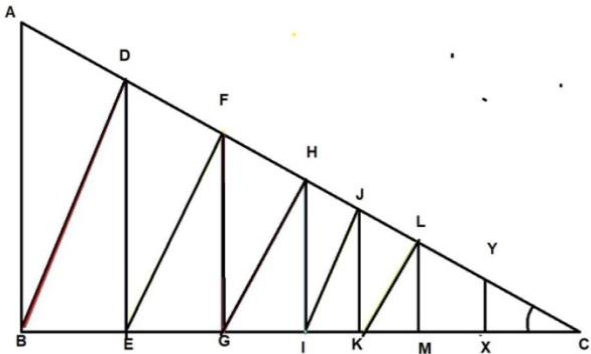

Figure 1 Computation of $x^n$, Euler's number e, m/n in equation $e^{m/n}$ and also reciprocal of a number

Considering the right angled $\triangle ABC$ and subsequent similar right angled triangles $BDC, DEC, EFC, FGC, \ldots$, so on, aforesaid perpendiculars can be written in the following way

$$BD = p_1 = AB \cos P \cos C,$$
$$DE = p_2 = AB \cos^2 C = P \cos^2 C,$$
$$EF = p_3 = AB \cos^3 C = P \cos^3 C, \quad (1)$$
$$FG = p_4 = AB \cos^4 C = P \cos^4 C,$$

…. so on and
$$XY = p_n = AB \cos^n C = P \cos^n C.$$

## 2. MATERIALS AND METHODS

In this section, computation of different operations using real numbers except zero or infinity as detailed in introduction, will be taken up one by one followed by fabrication and working of mechanical analogue calculator.

### 2.1. Computation of $x^n$

2.1.1. *When real number $x > 0$ and n is a positive integer but is neither zero nor infinity*

For computing $x^n$, when real number $x > 1$, $x$ can not be equated with $\cos C$ being more than 1. But when $1 > x > 0$, geometric figure can be drawn as $x$ can be equated with $\cos C$. When $x > 1$, I can write $x = 10^k X$, where $1 > X > 0$ and $k$ is any real number, positive or negative, then $X$ can be equated with $\cos C$ and geometric figure with angle $C$ corresponding to $\cos C = X$ can be constructed as explained in section 1. $x^n$ then equals to $10^{kn} X^n$. Perpendicular $BD$ from point $B$ is drawn upon hypotenuse $AC$ and another perpendicular $DE$ from point $D$, is drawn upon base $BC$. The process of drawing perpendiculars upon hypotenuse and, then upon base, is continued till total number of such perpendicular equals $n$ and length of $n^{th}$ perpendicular say $YX$ is then given by relation

$$YX/AB = p_n/P = \cos^n C = X^n \quad (2)$$

or

$$x^n = 10^{kn} X^n = 10^{kn} (p_n/P). \quad (3)$$

When length of perpendicular $AB = P = 1$, then $x^n = 10^{kn}(YX) = 10^{kn}(p_n)$. Value of product of $k$ and $n$, can be computed as will be described in section 2.2. If $x$ is negative and $n$ is odd integer, then $x^n$ will be negative and in all other cases, positive.

*Example*:

Compute $x^n = (32357)^{10}$. In this case, $(32357)^{10}$ can be written as $(0.32357)^{10} \cdot 10^{50}$, $\cos C$ is, then equated with $0.32357$ and triangle $ABC$ with angle $C$ corresponding to $\cos C = 0.32357$, is constructed. Thereafter, perpendiculars are drawn on hypotenuse $AC$ and base $BC$ as already described.

2.1.1.1. Alternate method of computation of $x^n$ when real number x > 1: In such cases, $1/x^n$ can be computed by equating $1/x$ with $\cos C$, geometric construction is same as Figure 1, but $x^n = 1/p_n$. Length $p_n$ can be measured and its reciprocal can be calculated as explained in section 2.2.2.1.2.

*Computation of $x^{-n}$, when real number x*

*is such that $0 < x < 1$ and $n$ is a positive integer but not zero*

$x^{-n}$ can be written as $(y)^n$, where $y = 1/x$. Computation of $y^n$ has already been explained in section 2.1.1. Value of $x^{-n}$ will be reciprocal of $y^n$. Computation of reciprocal of a real number has been explained in section 2.2.

*2.1.3. Computation of $x^{-n}$, when real number $x > 1$ but is not infinity and $n$ is a positive integer but not zero*

In such cases, $1/x$ can be equated with $\cos C$ and value of $(x)^{-n}$ is computed as explained in section 2.1.1.1. *In addition to this method, by writing $x^{-n} = (y)^n$ where $y = 1/x$, method of scaling up or down, can also be used as explained in section 2.1.1.*

*Examples:*

i) Compute $x^{-n} = (32357)^{-10}$. Quantity $(32357)^{-10}$ is written as $(1/32357)^{10}$ and $32357$ is written as $3.2357 \cdot 10^4$. Quantity $(1/32357)^{10}$ is written as $(1/3.2357)^{10} \cdot 10^{-40}$ and $\cos C$ is equated with $1/3.2357$. Rest procedure will be same as described for positive $n$.

*2.1.4. Algorithms or Pseudocode or Flowcharts*

With the purpose of forming a computable machine programmes and also for the clarification of the methods used, I will be describing the algorithms for different mathematical operations. Although these algorithms in the strictest sense, pseudocode or flow charts are additional to the contents of this paper, these do highlight the importance of this research paper in geometric computations of mathematical operations.

*2.1.5. Algorithm or Pseudocode or Flowchart for Computations of $x^n$, where $x$ and $n$ are real numbers positive or negative but neither zero nor infinity*

1. Given $x^n$. Check whether $n$ is positive or negative. If it is negative go to 3.
2. Write $x^n = 10^{kn}(X)^n$ so that $1 > X > 0$ and $k$ has an integer value positive or negative including 0. Then go to 4.
3. Write $n = -N$ and $x^{-N} = 10^{-kN}(X)^{-N}$ so that $1 > 1/X > 0$ and $k$ has an integer value positive or negative including 0. Then go to 5.
4. Construct triangle $ABC$ such that $AB = 1$, angle $C$ corresponds to $\cos C = X$ and angle $B$ equals to 90 degrees. From point $B$, draw perpendicular $BD$ on $AC$, from point $D$, draw perpendicular $DE$ on $BC$,..., so on. Number of these perpendiculars must be equal to $n$. Measure length of the $n^{th}$ perpendicular considering $BD$ as the first perpendicular. Let it be $XY$. Then go to 6.
5. Construct triangle $ABC$ such that $AB = 1$, angle $C$ corresponds to $\cos C = 1/X$ and angle $B$ equals to 90 degrees. From point $B$, draw perpendicular $BD$ on $AC$, from point $D$, draw perpendicular $DE$ on $BC$,..., so on. Number of these perpendiculars must be equal to $n$. Measure length of the $N^{th}$ perpendicular considering $BD$ as the first perpendicular. Let it be $XY$. Then go to 8,
6. Print result $XY(10^{kn})$ if $x$ is not -ve and $n$ is also not odd otherwise go to 7.
7. Print result $-XY(10^{kn})$ and go to 10.
8. Print result $XY(10^{-kN})$, if $x$ is not -ve and $N$ is also not odd otherwise go to 9.
9. Print result $-XY(10^{-kN})$.
10. Stop.

In proceeding sections, algorithms for each operation will not be explained for the sake of brevity but it can be formed following the procedure as detailed above in accordance with the specific need of the mathematical operation.

## 2.2. Operations of multiplication, division of two real numbers, geometric mean of two positive real numbers and also reciprocal of a number provided the numbers are neither zero nor infinity

*2.2.1. Reciprocal of a real number $p(10^k)$, when $0 < p < 1$*

Referring to Equations (1) and Figure 1, $BD = AB \cos C$ and given number is $p(10^{k_2})$, where $p > 1$. If $AB = 1$ and $\cos C = 1/p$, then $BD$ will equal $1/p$. Construct a triangle $ABC$ with angle $C$ corresponding to $\cos C = 1/p$, perpendicular $AB = 1$, base $BC$ and hypotenuse $AC$. Draw perpendicular $BD$ upon $AC$ meeting it at point $D$, then $1/p = BD(10^{-k_2})$.

*Example*:
Computation of reciprocal of electron charge. Electron charge is $-1.602176634 \times 10^{-19}$ Coulomb. Ignoring negative sign for time being and writing it as $1.602176634 \times 10^{-19}$, I construct triangle ABC with angle $C$ corresponding to $\cos C = 1/1.602176634$, perpendicular $AB = 1$, base $BC$, hypotenuse $AC$ and right angle at point $B$. From pint point $B$, I draw a perpendicular $BD$ upon $AC$ meeting it at point $D$, then $BD = 1/(1.602176634)$ and reciprocal of electron charge will be $-BD(10^{19})$ per Coulomb.

Another Method: Referring to Equations (1) and Figure 1,
$$BD/AB = \cos C$$
$$DE/DB = \cos C,$$
therefore, $BD^2 = (AB)(DE)$ or
$$BD = \sqrt{(AB)(DE)}. \quad (4)$$
Let $BD = 1, AB = P$ and $DE = p_2$, then $1 = Pp_2$ or
$$P = 1/p_2, \quad (5)$$
where $P > p_2$ and $0 < p_2 < 1$. From a point $D$, a perpendicular $DE$ of given length $p_2$ drawn upon a straight line $B'C'$ meeting it at point $E$. From the same point $D$, a straight line $DB$ of unit length is drawn meeting the straight line $B'C'$ at point $B$, thereafter, a straight line $A'C$ perpendicular upon $DB$, passing through point $D$ and meeting $B'C'$ at $C$ is drawn. A perpendicular $BA$ is raised from point $B$ meeting $A'C'$ at $A$, then $AB$ will measure $1/p_2$ reciprocal of $p_2(10^k)$ will equal $AB(10^{-k})$.

*2.2.2. Geometric mean of two positive (or negative) real numbers P and $p_2$, where $P > p_2$ and none of these is equal to zero*

To find geometric mean of two real positive (or negative) numbers $P$ and $p_2$, I will explain two methods. In first method, I will use trigonometry and geometry both and in second only geometry.

2.2.2.1. Method using both trigonometry and geometry or angle bisection method: Perpendicular $DE$ can be written equal to $(BD) \cos C$ and perpendicular $BD = (AB) \cos C$. That makes $DE = (AB) \cos^2 C$ and on substituting values of $DE$ and $AB$, I get $\cos^2 C = p_2/P$. Since according to cosine doubling formula, $\cos 2C = 2\cos^2 C - 1$, therefore,
$$\cos 2C = (2p_2 - P)/P. \quad (6)$$
Referring to Figure 2, if quantity $(2p_2 - P)$ happens to be negative, that means angle $2C$ is obtuse otherwise acute. Angle $XCY$ corresponding to $\cos 2C = (2p_2 - P)/P$ is drawn and is, then bisected to construct angle $XCY'$ which will equal angle $C$. A perpendicular $ED$ is raised upon straight line $XC$ intersecting straight line $CY'$ at point $D$ such that length of $ED = p_2$. From the same point $D$, another perpendicular $DB$ is raisedupon $CD$ to meet straight line $XC$ at point $B$. Thereafter, from point $B$, another perpendicular $BA$ is raised upon straight line $CB$ to meet at point A straight line $CY'$. Length $BD$ will then equal to $\sqrt{(AB)(DE)}$ or $\sqrt{Pp_2}$.

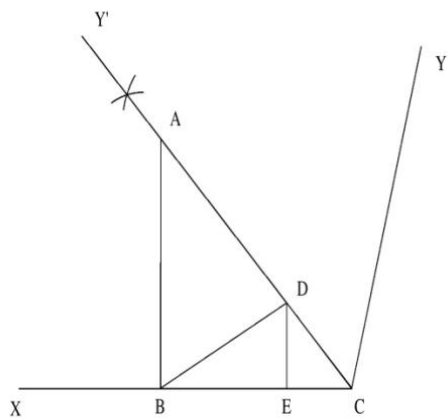

Figure 2 Trigonometric or angle bisection method for computing geometric mean

2.2.2.2. Method using geometry: Let $P = P_{k_1}(10^{k_1})$ and $p_2 = P_{k_2}(10^{k_2})$, where $P_{k_1} > P_{k_2}$, both $k_1$, $k_2$ are either both even or both odd integers and have values such that both $P_{k_1}$, $P_{k_2}$ have values between 0 and 1. In this method, referring to Figure 3, value of angle $C$ is explored so that length of perpendicular $DE$ equals given value $P_{k_2}$, also perpendicular $AB$ equals given value $P_{k_1}$ and $BD$ remains perpendicular upon hypotenuse $AC$.

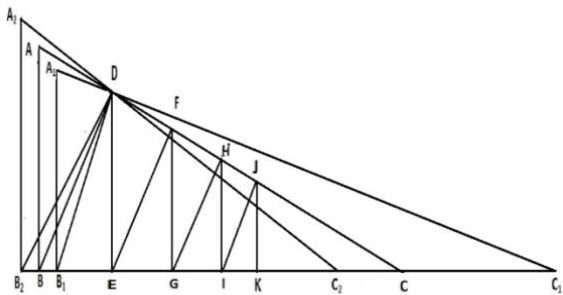

Figure 3 Operations of multiplication, division, geometric mean

As far as value of angle $C$ is concerned, in the first instance, it is assigned a value arbitrarily. Perpendicular $DE$ equal to value $P_{k_2}$ is drawn upon base $CB$, perpendicular $DB$ is drawn on hypotenuse $AC$ to meet at point $B$ the straight line $CE$, when extended and perpendicular $BA$ is drawn upon $CB$ at point $B$ to meet at point $A$ the straight line $CD$, when extended, then the length of perpendicular $BA$ is measured. Two cases will arise.

*First case,* when the length of perpendicular $AB$ happens to be more than the value $P_{k_1}$. Kindly refer to Figure 3. Length of $AB >$ $P_{k_1}$ indicates that the value of angle $C$ is more than the required value and thus needs reduction. *Hypotenuse $AC$ is slightly rotated with centre at $D$ to position $A_1C_1$ thus reducing the value of angle $BCA$ to angle $BC_1A_1$. From point $D$, perpendicular $DB_1$ is drawn upon hypotenuse $A_1C_1$ to meet base $BC_1$ at point $B_1$, from point $B_1$, perpendicular $B_1A_1$ is drawn upon $B_1C_1$ to meet titled hypotenuse $A_1C_1$ at $A_1$. If the length of perpendicular $B_1A_1$ is found equal to $P_{k_1}$, then geometric mean of $P$ and $p_2$ will be $B_1D\{10^{(k_1+k_2)/2}\}$.*

*Second case,* when the length of perpendicular $AB$ happens to be less than the value $P_{k_1}$. Kindly refer to Figure 3. That indicates the value of angle $C$ is less than the required value and needs enlargement. In that eventuality, *hypotenuse $AC$ is slightly rotated with centre at $D$ to position $A_2C_2$ thus increasing the value of angle $BCA$ to angle $BC_2A_2$. From point $D$, perpendicular $DB_2$ is drawn upon hypotenuse $A_2C_2$ to meet base $BC_2$ at point $B_2$ when the base is extended, from point $B_2$, perpendicular $B_2A_2$ is drawn upon $B_2C_2$ to meet titled hypotenuse $A_2C_2$ at point $A_2$. If the length of perpendicular $B_2A_2$ is found equal to $P_{k_1}$, then the length $B_2D$ is measured and the geometric mean of $P$ and $p_2$ will be $B_1D\{10^{(k_1+k_2)/2}\}$.*

2.2.3. *Product of two real positive numbers $P$ and $p_2$ where $P > p_2$ and neither of these is zero nor infinity*

The method described for determining the geometric mean is used for finding product of two numbers say $P$ and $p_2$. When geometric mean say $DB_1\{10^{(k_1+k_2)/2}\}$ is determined, then squaring this geometric mean will yield product of two numbers. Squaring of a real number has already been explained in section 2.1.

2.2.4. *Division of a real number by another real number*

Let $P = P_k(10^k)$ and $p = P_{k_2}(10^{k_2})$, where $1 > P_k > 0$ and $1 > 1/P_{k_2} > 0$.

Construct a triangle $ABC$ with angle $C$ corresponding to $\cos C = 1/P_{k_2}$ and perpendicular $AB = P_k$. From point $B$, draw a perpendicular $BD$ upon $AC$ meeting it at point $D$, then $P/p = BD(10^{k-k_2})$.

Second method: Let there be real numbers $P_k(10)^{k_1}$ and $p_{k_2}(10)^{k_2}$, where $k_1, k_2$ are real positive or negative integers such that $0 < P < 1$ and $P > p$. I want to calculate value of $p(10)^{k_2}$ divided by $P(10)^{k_1}$. Its value can be written as $(p/P)10^{k_2-k_1}$. Construction of the figure in this case will be different from Figure 1 and is explained hereinafter. A triangle $ABC$ with right angle at point $B$, base $BC = 1$ and perpendicular $AB = P$ is constructed. From a point $X$ on base $BC$, a perpendicular $XY = p$ is raised meeting hypotenuse $AC$ at point $Y$, then triangles $ABC$ and $YXC$ are similar therefore $CX/1 = p/P$ or length $CX$ will equal $p/P$ and $p(10)^{k_2}/P(10)^{k_1}$ will equal $CX(10)^{k_2-k_1}$.

## 2.3. Computation of value of $x^{1/n}$ where $\infty > x > 0$ and n is a real nu*mber*

Referring to Figure 1, $BD, EF, GH, IJ, KL$ ..., or $p_1, p_3, p_5, p_7, p_9, ...$, are $1^{st}, 3^{rd}, 5^{th}, 7^{th}, 9^{th}$ perpendicular upon hypotenuse $AC$ pertain to odd powers of $\cos C$ as their lengths are $(AB)\cos C$ $(AB)\cos^3 C$, $(AB)\cos^5 C$, $(AB)\cos^7 C$, $(AB)\cos^9 C$ respectively. Similarly, $DE, FG, HI, JK, LM, ...,$ or $p_2, p_4, p_6, p_8, p_{10}, ...$, are $2^{nd}, 4^{th}, 6^{th}, 8^{th}, 10^{th}...$ perpendicular upon base $BC$ pertain to even powers of $\cos C$ as their lengths are $(AB)\cos^2 C$, $(AB)\cos^4 C$, $(AB)\cos^6 C, (AB)\cos^8 C, (AB)\cos^{10} C$, ... respectively. Let $x^{1/n} = \cos C$, then $x = \cos^n C$. Let perpendicular say $XY = p_n = (AB)\cos^n C$. Assuming $x = AB, XY = 1$ and $x > 1$ then $1 = x \cos^n C$ or
$$(1/x)^{1/n} = \cos C. \quad (7)$$
That means, when $x > 1$, value of $1/\cos C$ requires determination. On the other hand, when $x < 1$, assuming $AB = 1$ and perpendicular $XY = x$, then
$$(x)^{1/n} = \cos C. \quad (8)$$

Since values of $x$ and $n$ are given, therefore, value of $\cos C$ requires determination for evaluating $x^{1/n}$. $x$ can also be scaled down or up by writing it as $(10^{kn}) X$ or $(10^{-kn}) X$ as the case may be, where $X < 1$. Depending upon odd, even parity of $n$, I will take up two cases.

### 2.3.1. When n is even integer and $\infty > x > 0$

First, the case, when $x > 1$ and $n$ is even, is being considered for computation of $x^{1/n}$. I assume $n = 4$ in the first instance and, then will generalise it to all positive even integers $n$ except $n = \infty$. Kindly refer to Figure 4. For construction, a horizontal straight line $BC$ is drawn and from a point $G$ on it, a perpendicular say $GF$ of unit length is raised upon it. Points $C$ and $F$ are joined and straight line $CF$ is extended, say to point $A$. Since $GF$ is perpendicular upon base $BC$, therefore, it will denote even power of $\cos C$. From point $F$, a perpendicular $FE$ is raised upon hypotenuse $AC$ intersecting base $BC$ at $E$. From point $E$, a perpendicular $ED$ is raised upon base $BC$ intersecting $AC$ at point $D$. From point $D$, a perpendicular $DB$ is raised upon hypotenuse $AC$ to intersect base $BC$ at point $B$ Then, from point $B$, a perpendicular $BA$ is raised upon base $BC$ intersecting hypotenuse at point $A$. Let $n = 4$, there will be four perpendicular namely $BD, DE, EF$ and FG in triangle $ABC$. Two cases will arise.

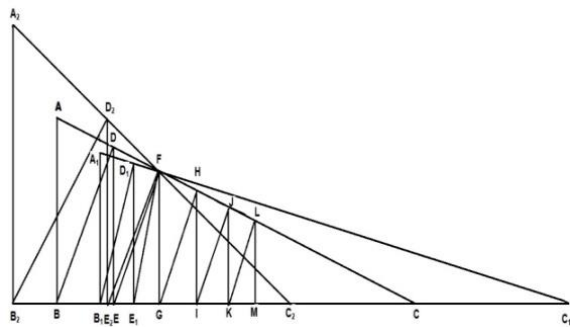

Figure 4 Computation of $x^{1/n}, x^{m/n}$, $ln(a)/ln(x)$ and $ln(x)/ln(a)$

a. If on measuring, length of perpendicular $AB$, it is found more than $x$, then *value of angle BCA needs to be reduced in such a way that hypotenuse AC passes through point F*. Kindly refer to Figure 4. Reduced angle C is now $BC_1F$ and hypotenuse is $A_1C_1$. Again perpendiculars are drawn upon base $BC_1$ and hypotenuse $A_1C_1$ and now if length $A_1B_1$ is found equal to $x$, then the value of $x^{1/4}$ is $A_1C_1/B_1C_1$.

b. If on measuring, length of perpendicular $AB$, is found less than $x$, then *value of angle BCA needs to be increased in such a way that hypotenuse AC passes through point F*. Enlarged angle C is now $BC_2F$ and hypotenuse is $A_2C_2$. Again perpendiculars are drawn upon base $BC_2$ and hypotenuse $A_2C_2$ and now if length $A_2B_2$ is found equal to $x$, then the value of $x^{1/4}$ is $A_2C_2/B_2C_2$. In the above case, $n$ was taken equal to 4. If $n$ has some other positive even integer value, then number of perpendiculars will be $n$ instead of 4 but the procedure will be same as has already been detailed.

Now the case, when $0 < x < 1$ and $n$ is even positive integer, is considered. Assuming $n$ equal to 4, geometric construction will be same as that of earlier case but the length of the perpendicular $FG$ will be $x$ and the length of perpendicular $AB$ will be unity. Rest steps will be same as that of earlier case but the result will be $\cos C_1$ or $\cos C_2$, that is, $x^{1/n}$ will equal $B_1C_1/A_1C_1$ or $B_2C_2/A_2C_2$ as the case may be.

When $x$ is of the form $l/m$ and $\infty > l/m > 1$ and $n$ (in this case it is 4) is even integer, then a straight line $BC$ is drawn and a point $G$ is taken upon it, from this point, a perpendicular say $GF$ of length $l$ is raised in stead of unit length. Rest procedure is same as described above except measurement of $AB$ or $A_1B_1$ or $A_2B_2$ as the case may be, should equal to value $m$ in stead of $x$.

When $x$ is of the form $0 < l/m < 1$ and $n$ (in this case it is 4) is even integer, then a straight line $BC$ is drawn and a point $G$ is taken upon it, from this point, a perpendicular say $GF$ of length $m$ is raised in stead of unit length. Rest procedure is same as described above except the measurement of $AB$ or $A_1B_1$ or $A_2B_2$ as the case may be, should equal to value $l$.

*If $x$ is negative and $n$ is even integer then $x^{1/n}$ is not computable for its solution in real numbers.*

2.3.2. *Computation of value of $x^{1/n}$ when $n$ is odd integer and $x$ is positive or negative real number but is neither 0 nor infinity*

Referring to Equations (1),
$DE = p_2 = p_1 \cos C = BD \cos C,$
$EF = p_3 = p_1 \cos^2 C = BD \cos^2 C,$
$FG = p_4 = p_1 \cos^3 C = BD \cos^3 C,$
…, so on,
$XY = p_{n+1} = p_1 \cos^n C = BD \cos^n C. (9)$

Here I am referring to perpendicular $BD$ of triangle $BDC$ in stead of perpendicular $AB$ of triangle $ABC$ that was being used earlier. Purpose of consideration of triangle $BDC$ and writing the equations in this fashion is to keep intact the similarity with the procedure explained for even value $n$. By referring to perpending $BD$ of triangle $BDC$, I will get odd numbered perpendiculars upon horizontal line $BC$ and will also be able rotate side $DC$ to alter angle $C$ for computation of $x^{1/n}$. If $x$ is negative quantity, positive $x$ will be used but the final result will be given negative sign.

Kindly refer to Figure 4. The steps for calculations of $x^{1/n}$ are same as explained earlier except the fact that perpendicular $BD$ will be used in place of perpendicular $AB$. Since perpendicular $AB$ of triangle $ABC$ is analogous to perpendicular $BD$ of triangle $BD$, perpendicular $FG$ in triangle $BDC$ will now pertain to $n = 3$. Kindly refer to Equations (9).

*Examples:* i) Computation of value of $(32157/121)^{1/16}$. Obviously $32157/121 > 1$ and can also be written as $32.157/0.121$. △ $ABC$ can be constructed by drawing a straight line $BC$ and raising a perpendicular $Y_8 X_8 = 0.121$ from any point $Y_8$ on base $BC$, drawing hypotenuse by joining points $C$ and $X_8$, then extending straight line $CX_8$. Procedure given in paragraph 2.3.1 is then followed. When perpendicular $AB$ or $A_1 B_1$ or $A_2 B_2$ will be equal to 32.157, then $(32.157/0.121)^{1/16} = AC/BC$ or $A_1 C_1 / B_1 C_1$ or $A_2 C_2 / B_2 C_2$.

ii) Computation of value of $(3.2157/121)^{1/4}$. $\triangle ABC$ can be constructed by drawing a straight line $BC$ and raising a perpendicular $Y_2 X_2 = 3.2157$ from any point $Y_2$ on base $BC$, drawing hypotenuse by joining points $C$ and $X_2$, then extending straight line $CX_2$. Procedure given in paragraph 2.3.1 is then followed. When perpendicular $AB$ or $A_1 B_1$ or $A_2 B_2$ will be equal to 121, then $(3.2157/121)^{1/4} = BC/AC$ or $B_1 C_1 / A_1 C_1$ or $B_2 C_2 / A_2 C_2$.

iii) Calculate $(-32157/121)^{-1/17}$. $\triangle ABC$ can be constructed by drawing a straight line $BC$ and raising a perpendicular $Y_9 X_9 = .121$ from any point $Y_9$ on base $BC$, drawing hypotenuse by joining points $C$ and $X_9$, then extending straight line $CX_9$. Procedure given in paragraph 2.4 is then followed. When perpendicular $B_1 X_1'$ or $B_1 X_1'$ or $B_2 X_2''$ will be equal to 32.157, then $-(32157/121)^{1/17} = -AC/BC$ or $-A_1 C_1 / B_1 C_1$ or $-A_2 C_2 / B_2 C_2$.

### 2.4. Computation of $x^{m/n}$

For such cases, where $m$ and $n$ are integers neither zero nor infinity and $\infty > x > 0$, $x^m$ is first computed by the method as described in section 2.1 and to avoid repetition, it is not described again. Let $x^m$ on geometric computation equals $y$, then $x^{m/n}$ will equal to $y^{1/n}$. For computing $y^{1/n}$, method has already been described in Section 2.3. Therefore, computing first $x^m$ and then $(x^m)^{1/n}$ will yield value of $x^{m/n}$. However, an alternate approach to compute $x^{m/n}$ can be adopted by writing it in the form $x^{m_1 + m_2/n}$ or $(x^{m_1})(x^{m_2/n})$ where $m = n \cdot m_1 + m_2$. Then $x^{m_1}$ can be computed, according to the method given in Section 2.1, and $x^{m_2/n}$ can be computed as already explained. Thereafter, product of $x^{m_1}$ and and $x^{m_2/n}$ is computed as explained in section 2.2. Easiness in computation by this method is attributed to the fact that lower value of $m_2$ will entail less operations. If exponent $m/n$ happens to be in decimals, then whole integer can be considered $m_1$ and decimal part converted to fraction of the form $m_2/n$. To avoid repetition, I am not describing algorithm for computation of $x^{m/n}$ which, in fact, is composite of algorithm of $x^m$ and $x^{1/n}$. If $-\infty < x < 0$, then $n$ must be an odd integer for compatibly of $x^{m/n}$.

### 2.5. Computation of integer $n$, when values of $x$ and $a$ in equation $x^n = a$ or $x^{1/n} = a$, are given

a. Assuming $n$ positive integer which is neither zero nor infinity, equation $x^n = a$, where $\infty > x > 1$ and $a$ is real number neither zero nor infinity, can be written as $(1/x)^n = 1/a$ so that $1/x$, being less than 1, can be equated to $\cos C$, where angle $C$ is enclosed by the sides $BC$ and $AC$ of triangle $ABC$, which has vertex at $A$, right angle at $B$ and perpendicular $AB$ of unit length. For geometric construction, refer to Figure 1. Length of each perpendicular, $AB, BD, DE, EF$ ... are compared one by one with $1/a$ till length of $N^{th}$ perpendicular say $XY$ equals to $(1/a)$, then $n = N$.

b. Assuming $n$ positive integer which is neither zero nor infinity, equation $x^n = a$, where $0 < x < 1$ and $a$ is real number neither zero nor infinity, can be written as $(x)^n = a$, then triangle $ABC$ is drawn with vertex at $A$, right angle at $B$, perpendicular $AB = 1$ and $\cos C = x$. Procedure as already described is used here till length of $N^{th}$ perpendicular say $XY$ equals to $(a)$, then $n = N$.

c. Assuming $n$ positive integer which is neither zero nor infinity, equation $x^{1/n} = a$, where $\infty > x > 1$ and $a$ is real number

neither zero nor infinity, can be written as $1/x = (1/a)^n$. Triangle $ABC$ is drawn with vertex at $A$, right angle at $B$, perpendicular $AB = 1$ and $1/a = \cos C$. The procedure is same as described in section 2.5.a for $x^n = a$ except perpendiculars $BD, DE, EF\ldots$ will be compared with $1/a$.

d. If $n$ is positive odd integer which is neither zero nor infinity, equation $x^{1/n} = a$, where $-\infty < x < 0$ and $-\infty < a < 0$, then assuming $y = -x$, and $b = -a$ will yield $(y)^n = b$. Now $\infty > y > 0$ such cases have already been taken up.

**2.6. Computation of $m/n$, when values of $x$ and $a$ in equation $x^{m/n} = a$ are known and $m, n, x, a$ are positive numbers such that none of these is zero or infinity.**

When $1 < x < a$ in equation $x^{m/n} = a$, then obviously $m > n$ and can be written as $m/n = m_1 + n_1/n$. Equation $x^{m/n} = a$ takes the form

$$(1/x)^{m/n} = 1/a. \quad (10)$$

Construction of geometric figure for this computation is same as Figure 1 but triangle $ABC$ has angle $C$ corresponding to $\cos C = (1/x)$. Length of each perpendicular $AB, BD, DE, EF \ldots$ is compared one by one with $1/a$. Let the length of $(N+1)^{th}$ perpendicular designated as $p_{N+1}$ be less than $1/a$ and the length of $(N)^{th}$ perpendicular $p_N$ be more than $1/a$, then $m_1 = N$ and the ratio of $1/a$ and length of $N^{th}$ perpendicular $p_N$ equals to $(1/x)^{n_1/n}$, where $n > n_1$ or $1/(a\, p_N) = (1/x)^{n_1/n}$.

Therefore,

$$1/(a \cdot p_N)^{n/n_1} = (1/x). \quad (11)$$

Let $n/n_1 = m_2 + n_2/n_1$. Since value of $a$ is known, length $p_N$ can be measured, value of $x$ is given, therefore, Equation (11) is similar to Equation(10).. Geometric Figure 1 is reconstructed. Now triangle $ABC$ with vertex at $A$, right angle at $B$, $AB = 1$ has angle $C$ corresponding to $\cos C = 1/(a\, p_N)$. Fresh perpendicular $BD, DE, EF, FG \ldots$ are drawn. Again length of each perpendicular $AB, BD, DE, EF \ldots$ is compared one by one now with $1/x$. Let the length of $(N_1 + 1)th$ perpendicular designated as $p_{N_1+1}$ be less than $1/x$ and the length of $(N_1)th$ perpendicular designated $p_{N_1}$ be more than $1/x$, then $m_2 = N_1$ and the ratio of $1/x$ and the length of $(N_1)th$ perpendicular $p_{N_1}$ equals to $1/(a \cdot p_{N_1})^{n_2/n_1}$ where where $n_1 > n_2$. That is $1/(x \cdot p_{N_1}) = 1/(a \cdot p_N)^{n_2/n_1}$ or

$$1/(x \cdot p_{N_1})^{n_1/n_2} = 1/(a\, p_n). \quad (12)$$

Since values of $a, x$ are known, lengths $p_N$, $p_{N_1}$ can be measured, therefore, Equation (12) is similar to Equation (10) and can be processed for calculating value of $N_2$ in the same way as value of $N_1$ was determined. In this way, values of $N_3, N_4, N_5, \ldots$ can be determined and $m/n$ can be expressed by the following equation

$$\frac{m}{n} = N + \cfrac{1}{N_1 + \cfrac{1}{N_2 + \cfrac{1}{N_3 + \cfrac{1}{N_4 + \cdots}}}}. \quad (13)$$

Equation (13) is expression of $m/n$ in *continued fractions. This continued fraction will be of fixed terms if $m/n$ is a rational quantity and of infinite terms if $m/n$ is irrational quantity.* For example, if $m/n$ which is to be found out, is a rational quantity say $1971/181$, then $N = 10, N_1 = 1, N_2 = 8, N_3 = 20$ and $(1971/181) = 10 + \cfrac{1}{1+\cfrac{1}{8+\cfrac{1}{20}}}$. It may please be noted that $N, N_1, N_2, \ldots$ can have any positive integer value including 0 except infinity. It is also evident from the example that continued fraction has fixed terms. If $m/n$ is an irrational quantity and as irrational quantity can be at the most approximated with some value but can be represented by continued fractions with infinite terms, for example,

$$\sqrt{2} = 1 + \cfrac{1}{2+\cfrac{1}{2+\cfrac{1}{2+\cfrac{1}{2+\cdots}}}}.$$

Therefore, more the number of terms, more accurate the value of $m/n$ would be. Also it is explicit from this example that the fraction is continuing indefinitely, when $m/n$ is irrational. Once the continued fraction is determined, its value can be computed as it involves only divisions and computation of division of two numbers have already been explained.

*2.6.1. Computation of $m/n$ when equation $x^{m/n} = a$, has $x > 1$, and also $x > a$*

When the equation $x^{m/n} = a$, has $x > 1$, and also $x > a$, then obviously $n > m$ and $n/m$ can be expressed as $n_1 + m_1/m$. Equation $x^{m/n} = a$, then takes the form

$$(1/a)^{n/m} = 1/x, \qquad (14)$$

where $\cos C = 1/a$. Equation (14) is same as Equation (10), therefore, procedure applied to Equation (10) will also be applicable to Equation (14).

*2.6.2. Computation of $m/n$ when in equation $x^{m/n} = a, a < x < 1$*

When in equation $x^{m/n} = a, a < x < 1$, then obviously $m > n$ and $m/n$ can be expressed as $m/n = m_1 + n_1/n$. Then the equation, $x^{m/n} = a$, is same as Equation (10), therefore, procedure applied to Equation (10) will also be applicable to this equation.

*2.6.3. Computation of m/n when equation $x^{m/n} = a$, has $x < 1$ and also $x < a$*

When equation $x^{m/n} = a$, has $x < 1$ and also $x < a$, then obviously $n > m$ and $n/m$ can be expressed as $n_1 + m_1/m$. Equation $x^{m/n} = a$ takes the form $a^{n/m} = x$, where $\cos C = a$. This equation is same as Equation (10, therefore, procedure applied to Equation (210) will also be applicable to this equation.

*2.6.4. Computation of $m/n$ in equation $x^{m/n} = a$, when values of $x$ and $a$ are too large to construct a geometric figure and $m, n, x, a$ are positive numbers but are neither zero or infinity*

When values of $x$ and $a$ too large to construct a geometric figure due to the fact that $\cos C$ either approaches 0 or 1, the method detailed in section 2.7 becomes difficult to apply. In such a situation, change of base is resorted to and $x^{m/n} = a$ is written $m/n = (\ln a)/(\ln x)$. Assuming $p = \ln a$ and $q = \ln x$, these equations can be written $a = e^p$, $x = e^q$ and $m/n = p/q$. Since Euler's number for practical purposes is taken as 2.7182818 and is more than 1, therefore, Equation $a = e^p$, and $x = e^q$ are transformed to $1/a = (1/e)^p$ and $1/x = (1/e)^q$. Now triangle $ABC$ with angle $C$ corresponding to $\cos C$ is always constructible and value of $p, q$ in equations, $1/a = (1/e)^p, 1/x = (1/e)^q$ can be determined following the procedure explained in section 2.6. Let value of $p$ and $q$ on computation be

$$p = P + \frac{1}{P_1 + \frac{1}{P_2 + \cdots}} \quad \text{and} \quad q = Q + \frac{1}{Q_1 + \frac{1}{Q_2 + \cdots}},$$

then

$$\frac{m}{n} = \frac{p}{q} = \frac{P + \frac{1}{P_1 + \frac{1}{P_2 + \cdots}}}{Q + \frac{1}{Q_1 + \frac{1}{Q_2 + \cdots}}}. \qquad (15)$$

## 2.7. Approximation of Value of Euler's Number 'e'

Euler's number $e$, by definition is $\lim_{n \to \infty}(1 + 1/n)^n$ (Boyer and Merzbach 1991). Since it is greater than 1, therefore, equation for Euler's number can also be written as $1/e = \lim_{n \to \infty} 1/(1 + 1/n)^n$. Kindly refer to Figure 1. If perpendicular $AB = 1$ and angle $C$ corresponds to

$$\cos C = 1/(1 + 1/n),^n \qquad (16)$$

where $1/n$ is very small quantity tending to zero, then angle $C$ will also tend to zero. But for practical purposes, angle $C$ is taken having some value though extremely small. In that case, reciprocal of length of n$^{th}$ perpendicular will be given by equation

$$1/p_n = 1/\cos^n C = \lim_{n \to \infty}(1 + 1/n)^n = e.$$
$$(17)$$

## 2.8. Computation of natural logarithm and antilogarithm of a real number

If natural logarithm of a quantity say $a$ equals to $n$, that is $\ln a = n$, then $n$ is given by relation $e^n = a$. In section 2.5, detailed procedure has been given for such computation, where $x$ will have to be replaced by $e$. The procedure is not reproduced here to avoid repetition. Computation of antilogarithm of quantity $n$ can be defined as that value of $x$ such that $e^n = x$. In section 2.1, method of computing $x^n$ is already given. This method will also be applicable here, if $x$ is replaced with $e$. Kindly refer to sections 2.1, 2.3 and 2.4.

**2.9. Mechanical Analogue Calculator**

With the purpose of explaining how scientific calculations can be performed by simple model of mechanical analogue calculator using telescopic tubes, the details are given below. Further, it has already been explained that by this method of geometric computations, the calculations are error free but due to simplicity of this model, the computation may be approximate. That also does not mean that mechanical analogue calculator has intrinsic error, it can always be made as a state of art instrument but that requires more of engineering skill than mathematical one whereas this paper is mathematics oriented.

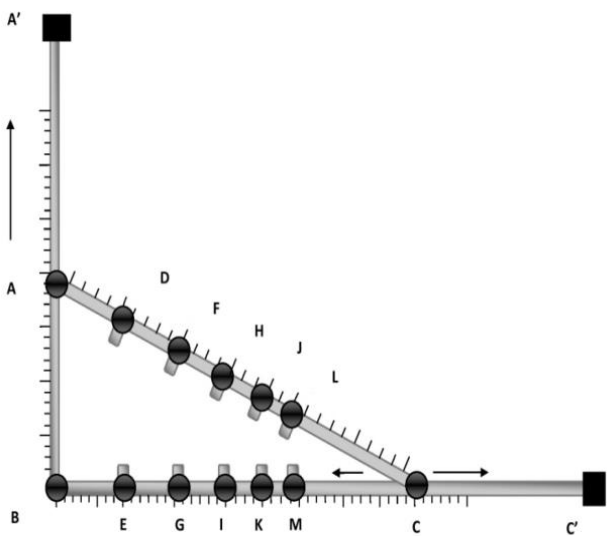

Figure 5a Mechanical Analogue Calculator in collapsed position

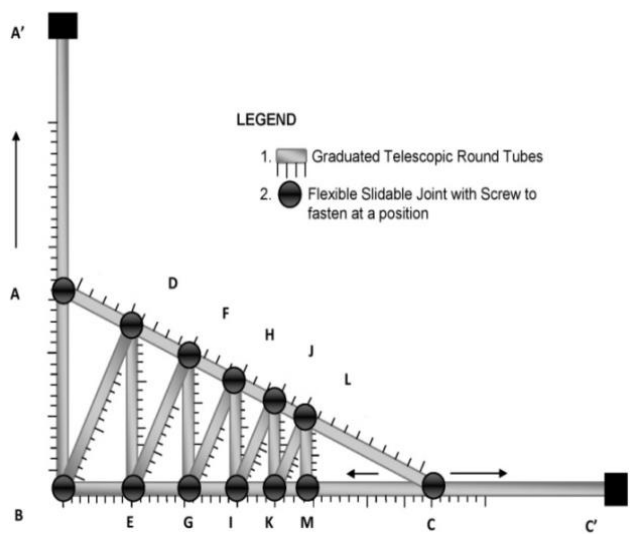

Figure 5b Mechanical Analogue Calculator in operating position

*2.9.1. Parts of mechanical analogue computing device*

Based on the theory of geometric computation using right angled triangles as described in the foregoing sections, a mechanical device has been presented for carrying out computations. It is submitted that French mathematician *Rene Descartes* while explain drawing of logarithmic curves, used somewhat this type of device (Smith and Latham 1954) [14]. This device has three main arms $A'B, BC'$ and $AC$, however it has $n$ perpendicular arms $B'D$, $D'E$, $E'F$, $F'G, G'H, H'I, I'J, J'K, K'L, L'M$ ... where value of $n$ denotes the maximum exponent that will be used in calculations. In total, the mechanical calculator has $n + 3$ arms, each arm has sufficient length say two metres or more. Longer the length, more accurate the results would be. Out of three main arms, two arms $A'B$ and $BC'$ are connected perpendicularly at joint $B$ which is moveable on arm $BC'$. Kindly refer to Figure 5a. The device shown in Figure 5a, has perpendicular arms $n = 10$ meaning thereby that its computability is limited to exponent 10. The arms are in shown collapsed position in Figure 5a, whereas these are in operating condition in Figure 5b.

*2.9.2. Structure and assembly*

2.9.2.1. Arms: Arms used in this device are telescopic graduated metal tubes, which fit into one another and are thus extendable and also collapsible, when pulled out or pushed in. These are designed so that their lengths can be read on the scale marked on them. This feature provides adjustability and readability for their lengths, which are essentially required for operability of the device for computations. Perpendicular arms $B'D$, $E'F$, $G'H, I'J, K'L$ have their fixed ends at $D, F, H, J, L$ respectively attached with joints $D, F, H, J, L$ slidable on hypotenuse A$C$. Other ends $B'$, $E'$, $G', I', K'$ of these arms are extendable to points $B, E, G, I, K$ respectively, where they are push fitted in joints $B$, joint $E$, joint $G$, joint $I$, joint $K$ respectively. These joints have these arms fixed in such a way that these are always perpendicular upon hypotenuse $AC$. I give explanation for notation that will be used in the text. *Apostrophe sign used after an alphabet means, that end of the arm is adjustable by compression or expansion, it also denotes that the end will finally be extended or compressed to same alphabet without apostrophe sign. For example, end A' will be expanded to joint A.*

Similarly, perpendicular arms $D'E, F'G$, $H'I, J'K, L'M$ have their fixed ends at $E$, $G$, $I, K, M$ respectively attached with joints $E$, $G, I, K, M$ slidable on base $BC$. Other ends $D'$, $F'$, $H', J', L'$ of these arms are extendable to points $D, F, H, J, L$ respectively, where these are push fitted in joints $D$, joint $F$, joint $H$, joint $J$, joint $L$ respectively.

2.9.2.2. Joints: The joints as stated above are slidable effortlessly on base $BC'$, hypotenuse AC and perpendicular $A'B$. In addition to fixation of perpendicular tubes upon these, these joints are designed to have mechanism so that the ends of a arms like $B', D', E'$ ... can be push fitted in them. After the arms had been adjusted for mathematical operation, the joints can be screwed to fasten them at the set positions so as to avoid unintentional disturbance.

2.9.2.3. Working: Initially all the arms except arm $A'B$ perpendicular upon $BC'$ and hypotenuse arm $AC$, are in collapsed position as shown in Figure 5a. For assembling triangle $ABC$ with angle $C$ corresponding to $cos\ C = x$, the length of arm $AC$ is adjusted to unity and the length of base $BC$ to $x$. Kindly refer to Figure 5b. End $A$ of hypotenuse arm $AC$ is detached from the joint $A$ and joint $A$ slid on $BA$ so that $AB = 1$. Joint $C$ with arm CA is slid on $BC'$, arm $CA$ is expanded or compressed *keeping the value of angle $C$ intact* and its end $A$ is push fitted in joint $A$. Joints $A$ and $C$ are then fastened with screws. Now angle $C$ corresponds to $cos\ C = x$ and length of $AB = 1$. Length of the arm $DB'$ is extended and the joint $D$ is slid upon hypotenuse $AC$ so that the extended end B' could be push fitted in joint $B$ and, then joint $E$ is slid on base $BC$ and the end $D'$ is extended so that end $D'$ is push fitted in joint $D$. In this way, all the perpendiculars are fitted to assemble analog calculator for performing computation of $x^n$. All the joints are screwed to fasten them at their set positions. Reading of lengths of perpendicular tubes will be noted. Value of $x$, $x^2$, $x^3$, $x^4$, $x^5$, $x^6$, $x^7$, $x^8$, $x^9$ or $x^{10}$ will equal to the length of $BD$, $DE$, $EF$, $FG$, $GH$, $HI$, $IJ$, $JK$, $KL$ or $LM$ respectively. I will take up examples here-under to explain how the calculator is operated for performing different mathematical operations. I submit, before performing a computation, all the perpendicular arms are set in collapsed position except arm $A'B$ perpendicular upon $BC'$. Kindly refer to Figure 5a. If this is not mentioned while explaining the operation, it would be assumed that all the perpendicular arms are in collapsed positions.

*2.9.3. Computation Using Mechanical Analogue Calculator*

2.9.3.1. Computation of $x^n$: Computation of electron charge $-1.602176634 \times 10^{-19}$ raised to power 7 will be performed using the calculator. Ignoring negative sign of the charge for time being, I write electron charge as $0.1602 \times 10^{-18}$, since numerals going beyond 4 digits place after decimal point requires big sized calculator, I take the value upto four decimal points. $x^n$, then equals to $(.1602)^{17}(10)^{-126}$. How multiplication of 18 with 7 is done to obtain 126, is explained in section 2.2. For assembling the device so that angle $C$ corresponds to $\cos C = .1602$, joint $C$ is slid upon $BC'$ to a position where $BC = 0.1602$. Length of $AC$ is made equal to 1.000 by compressing or expanding $AC$. Keeping intact the value of angle C, end A of the arm AC is detached and the length of AB is adjusted to unity. Perpendicular arms $BD, DE, EF, \ldots$ are adjusted and fitted as explained in section titled *Working*. The length of $7^{th}$ perpendicular $JK$ counting from $BD$ as first perpendicular is noted, then electron charge raised to power 7 is $-10^{-126}(JK)$.

2.9.3.2. Geometric mean and product of two numbers: Geometric mean of mass of earth $0.5972 (10^{25})$ Kg and mass of moon $0.7348(10^{23})$ Kg will be computed. At joint $E$, end $D'$ of arm $ED'$ is expanded and end $D'$ is push fitted in joint $D$ in such a way that $ED = 0.5972$. Joint $A$ is slid upon $A'B$ so that $AB = .7348$. Joint $E$ with arm $ED$ is slid on $BC$ and length $DB'$ is adjusted so that end $B'$ touches joint $B$. End $B'$ is push fitted in joint $B$. At joint $D$, arm $AC$ is expanded and end $A$ of arm $AC$ is push fitted in joint $A$ and end $C$ is push fitted in joint $C$. Length of $BD$ will be read and then geometric mean and product of of mass of earth and mass of moon is $BD\{10^{(25+23)/2}\}$ or $BD(10^{24})$ and $BD^2(10^{48})$ respectively.

2.9.3.3. Division of two numbers: For computing division of mass of earth by mass of moon, angle $C$ will be assembled corresponding to $\cos C = 1/7.348$. Join $C$ will be slid on $BC'$ so that $BC = 1$, joint $C$ is fixed at this position and joint $A$ is slid on $BA'$ while expanding $AC$ so that $AC = 7.348$. End $A$ of hypotenuse $AC$ is detached from the joint $A$. End $C$ of arm $AC$ is fastened to joint $C$ so that angle C is not disturbed. Joint $A$ is slid on $BA'$ so that $AB = .5972$. Joint $C$ is slid on $BC'$ keeping intact the value of angle $C$ so that end $A$ of the expanded $AC$ touches joint A where it is push fitted. At joint $D$, extend end $B'$ while sliding joint $D$ so that extended end $B'$ touches joint $B$, where it is push fitted in joint $B$. Length of $BD$ is read, then division of mass of earth by moon is $BD(10^{25-22})$ or $BD(10^3)$.

2.9.3.4. Computation of $x^{1/n}$: For computing $\{0.5972 (10^{25})\}^{1/6}$ where $0.5972 (10^{25})$ is the mass of the earth, I write it as $(5.972)^{1/6}(10^4)$ and consider $AB = 5.972$ and sixth perpendicular $(p_6), HI = 1.000$.. Procedure given in section, *Working*, is followed and length of $AB$ is fixed equal to 5.972. Length of $HI$ is read, if it is not equal to one, then keeping length $AB = 5.972$, value of angle $C$ is changed and the process is repeated till $HI$ is found equal to 1, then the lengths of $AC$ and $BC$ are read and computation of $\{0.5972 (10^{25})\}^{1/6}$ will equal to $10^4(AC/BC)$.

2.9.3.5. Computation of $x^{m/n}$: For computing $\{0.5972 (10^{25})\}^{19/7}$, I write it as $(.5972)^{19/7}(10)^{-1/7}(10^{68})$. Term $(.5972)^{19/7}$ and $(10)^{-1/7}$ are computable and their product multiplied by $(10^{68})$will give the result.

2.9.3.6. Computation of *m/n* when $x$ and $a$ are given in the equation $x^{m/n} = a$:Let us compute $m/n$ in equation $(98)^{m/n} = 151$ Since both 98 and 151 are too large quantities to draw angle $C$ corresponding to their reciprocals, therefore, method of $m/n = p/q = \ln(151)/\ln(98)$ will be resorted to and equation $(1/e)^p = 1/151$ and $(1/e)^q = 1/98$ will be considered for computation of $p$ and $q$. In mechanical calculator, joint $C$ and joint $B$ will be adjusted that length $CB = 1$ and $AC = e$ and perpendicular arms are constructed as detailed in section *Working*. Let length $1/151$ is less than $P^{th}$ perpendicular $P_p$ and more than $(P+1)^{th}$ perpendicular $P_{p+1}$, then value of $P$ is noted down. It gives new equation, $1/(151P_p)^{n/n_1} = 1/e$ for computation of $n/n_1$. Now angle $C$ will be constructed corresponding to $\cos C = 1/(151P_p)^{n/n_1}$ and value of $P_1$ in equation $1/(151P_p)^{n/n_1} = 1/e$, can be computed. Similarly, $P_2, P_3, \ldots$ can also be computed following the same procedure. In this way, value of $p = P + \frac{1}{P_1 + \frac{1}{P_2 + \cdots}}$, will be determined. Following the same procedure, value of $q = Q + \frac{1}{Q_1 + \frac{1}{Q_2 + \cdots}}$ is also determined. From these determinations, value of $m/n = p/q$ is finally worked out. It is observed that some computations involve a large number of steps, although laborious, it does not undermine the flawless computability of the device. On the contrary, it stresses the need of application of better mechanical engineering for improving its operability. Also, the mechanical computer displayed in Figure 5a and 5b, is the simplest form of mechanical computing machine or one can say, model for accomplishing the task of complicated mathematical operations. There is no doubt that with application of advanced mechanical engineering, improved form of the device can be put forth.

It is also clear that it is the distance measurement, after the device is set for operation, that matters. As far as distance measurement of perpendicular arms are considered, I have used graduated telescopic tubes. Evidently, the error in length measurement will culminate in error in results. This aspect I have explained in section titled *Accuracy and Precision.* No doubt the use of mechanical telescopic tubes, can also be done away with but that needs use of proximity sensors or any other electronic distance measuring instrument in place of telescopic tubes. This mechanical machine is manually operated and there is a need for automation of the device by servo control system that will also offset the error but this requires electronic Instrumentation. Incorporation of such salient features would make the machine self operating but that would be more of *Electronic Instrumentation Engineering* than *Geometric Computation.* What I want to emphasise, is that incorporating instrumentation engineering would do value addition to the mechanical calculating machine. But mechanical machine has its own indispensability on account of the fact that it does not require external power source (battery).

**2.10. Accuracy and Precision**

Length measuring instruments are *Ruler, Vernier Calliper, Micrometer, Electron Microscope, Transmission Electron Microscope.* A vernier caliper can measure a distance up to 0.01 mm; a bench micrometer up to 0.5 micrometer; a light microscope up to 0.2 micrometer, a Transmission Electron Microscope up to 0.1 nanometer. These instruments, on account of their resolving power, impose limitations on length measurements affecting accuracy of calculations.

**3. RESULTS AND DISCUSSION**

In a right angled triangle $ABC$ with right angle at point $B$, vertex at point $A$, perpendicular as $AB$ equal to $P$, base as $BC$ and hypotenuse as $AC$, a perpendicular of length say $p_n$ drawn upon base or hypotenuse from a point in hypotenuse or base as the case may be, the length of the perpendicular $p_n$ is given by relation,
$$p_n = P \cdot \cos^n C, \qquad (18)$$
where angle $C$ is angle $BCA$, $n = N_1 + n_1/n_2$, $N_1$, $n_1$ are integers including zero except infinity and $n_2$ is also an integer but not zero. Also if from point $B$, a perpendicular $BD = p_1$ is drawn upon hypotenuse meeting it at point $D$, then $p_1 = P \cdot \cos C$ and if from point $D$, a perpendicular $DE = p_2$ is drawn upon base $BC$ meeting it at point $E$, then $p_2 = P \cdot$

$cos^2 C$. In this way, if the practice of drawing perpendiculars is continued, then length of n[th] perpendicular $p_n$ will be given by relation, $p_n = P \cdot cos^n C$. But here, unlike $n = N_1 + n_2/n_1$, $n$ is an integer and thus it is a special case, where $n_2 = 0$. From these equations, we obtain $p_1^2 = P^2 \cdot cos^2 C = P(P \cdot cos^2 C) = P \cdot p_2$. Or $p_1 = \sqrt{P \cdot p_2}$, that is $p_1$ is geometric mean of $P$ and $p_2$ and also $p_1^2$ is product of $P$ and $p_2$. If length of $p_1$ is fixed as 1, then $P = 1/p_2$. Also since $p_n/P = cos^n C$ and if value of $P$ is fixed as 1, then computed value of $cos^n C$ equals to length of n[th] perpendicular. In this way, geometric mean, product of two numbers, reciprocal of a number and n[th] power of a quantity can be computed provided quantity whose n[th] power is to be computed is scaled down to less than 1 on account of the fact, $cos C \leq 1$. If $cos C = 1/(1 + 1/n)^n$ where $1/n$ is very very small quantity, then $cos^n C = 1/e = \lim_{n \to \infty} 1/(1 + 1/n)^n$. But $p_n = P cos^n C$ and if value of $P$ is taken as 1, then $p_n = cos^n C$ or $p_n = cos^n C = 1/e$. In this way, value of Euler's number can be computed. Also n[th] root of a quantity written in the form of $(p_n/P)^{1/n} = cos C$ can be computed by varying value of angle $C$ so that it should satisfy $p_n/P = cos^n C$.

For computing value of $n$ in equation, $p/P = cos^n C$, where $n = N_1 + n_2/n_1$, length of perpendicular $p_{N_1}$ corresponding to integer $N_1$ is found from condition when $p_{N_1} > P cos^n C$ and $p_{N_1+1} < P cos^n C$. At that stage, $p/p_{N_1} = cos^{n_2/n_1} C$ and by writing this equation as $cos C = (p/p_{N_1})^{n_1/n_2}$ where values of $cos C$, $p$ are given and value of $p_{N_1}$ can be measured from geometrical drawing. Obviously, $n_2 < n_1$, therefore, $n_1/n_2$ can be written as $N_2 + n_3/n_2$ where $N_2, n_3$ and $n_2$ are integers. Repeating the steps described for determining the value of $N_1$, value of $N_2$ can also be computed. In this way, value of $n$ can be computed from the continued fraction $N_1 + \frac{1}{N_2+\frac{1}{N_3+\cdots}}$. Whether the fraction continues finitely or infinitely depends upon the value of $n$. Method for computing reciprocal of a quantity is already known, therefore, continued fractions can be evaluated.

## 4. CONCLUSIONS

If in a triangle $ABC$ with base $BC$, perpendicular $AB$, hypotenuse $AC$ and right angle at point $B$, perpendiculars $BD, EF, GH$ are drawn upon hypotenuse $AC$, perpendiculars $DE, FH, HE, ...$ are drawn upon base $BC$, then $p_1 = BD = P cos C$, $p_2 = DE = P cos^2 C$, $p_3 = EF = P cos^3 C, ..., p_n = X_n Y_n = P cos^n C$.
Value of $x^n$ is the length of perpendicular $p_n$, when $0 < x < 1$ and $AB = 1$.
Value of $x^{1/n}$ is $1/cos C$ or $x/p_1$, when $x = AB$, $p_n = 1$ and $x > 1$.
By combining above stated both operations, value of $x^{m/n}$ is calculated.
Also since $p_1 = (P p_2)^{1/2}$, geometric mean of integer $m$ and $n$ is length $BD$, when $AB = m$, $DE = n$ and $m > n$.
When $1 < x < a$ and values of $x$ and $a$ are given in equation $x^{\frac{m}{n}} = a$, then $\frac{m}{n} = N_1 + \frac{1}{N_2+\frac{1}{N_3+\cdots}}$, when $AB = 1$, $cos C = (1/x)$, $p_{N+1} < 1/a$, $p_N > 1/a$, in subsequent triangle, $A'B' = 1$, $cos C' = 1/(a p_N)$, $p_{N_1+1} < 1/x$, $p_{N_1} > 1/x$ so on. If $m/n$ is irrational, then fraction will be continuous otherwise it will halt after fixed term. This geometric method can also be utilised for solving polynomial equations for their real roots having values between $-1$ and $+1$. Real roots having values more than one of the equation $a_0 x^n + a_1 x^{n-1} + a_2 x_{n-2} + \cdots + a_n = 0$, where coefficient $a_0, a_1, a_2, ..., a_n$ have fixed real values and $x$ is a variable, can be scaled down by substituting $x = 1/y$ and that makes applicable this geometric method to solving such equations (Wadhawan 2023).